\documentclass[11pt,a4paper]{article}
\usepackage[utf8]{inputenc}
\usepackage{amsmath,amssymb,amsthm,mathtools}
\usepackage{graphics} 
\usepackage{graphicx}  
\usepackage{epstopdf}
\usepackage[colorlinks=true]{hyperref}
\usepackage{tikz}
\usepackage{makeidx}
\usepackage{subfigure}
\usepackage{caption}
\usepackage{setspace}
\usepackage{a4wide}
\usepackage{bbm}
\graphicspath{{figures/}}
\theoremstyle{plain}
\newtheorem{theorem}{Theorem}[section]

\theoremstyle{definition}
\newtheorem{definition}[theorem]{Definition}

\theoremstyle{remark}

\newcommand{\R}{\mathbb{R}}
\newcommand{\G}{\mathcal{G}}

\newcommand{\vecsym}[1]{\boldsymbol{#1}}
\renewcommand{\vec}[1]{\mathbf{#1}}
\newcommand\norm[1]{\left\lVert#1\right\rVert}
\newcommand{\half}{\frac{1}{2}}
\renewcommand{\intercal}{T}

\begin{document}
 
\title{Recent Trends on Nonlinear Filtering for Inverse Problems}

\author{Michael Herty \footnotemark[1], \; Elisa Iacomini\footnotemark[1], \; Giuseppe Visconti\footnotemark[4]}

\footnotetext[1]{Institute for Geometry and Applied Mathematics, RWTH Aachen University,
	Templergraben 55, 52064 Aachen, Germany (herty@igpm.rwth-aachen.de, iacomini@igpm.rwth-aachen.de)}
\footnotetext[4]{Department of Mathematics, Sapienza University of Rome, P.le Aldo Moro 5, 00185 Rome, Italy (giuseppe.visconti@uniroma1.it)}

\date{}
\maketitle

\begin{abstract}
\noindent
Among the class of nonlinear particle filtering methods, the Ensemble Kalman Filter (EnKF) has gained recent attention for its use in solving  inverse problems. We review the original  method and discuss recent developments in particular in view of the limit for infinitely particles and extensions towards stability analysis and multi--objective optimization. We illustrate the performance of the method by using test inverse problems from the literature.
\end{abstract}



\section{Introduction} \label{sec:intro}

This review paper focuses on the Ensemble Kalman Filter applied to general inverse problems. In this context, some literature also uses the term  Ensemble Kalman Inversion (EKI). The method itself belongs to the class of particle methods and it is an iterative method for solving inverse problems. The method was originally introduced in~\cite{iglesiaslawstuart2013} for unconstrained minimization problems, and recently extended also to the presence of different types of constraints~\cite{ChadaShillingsWeissmann2019,HertyVisconti2020,Stuart2019CEnKF}. The original EnKF has already been introduced more than ten years ago~\cite{BergemannReich,ChenOliver,EmerickReynolds,Evensen1994} as a discrete time method to estimate state variables and parameters of stochastic dynamical systems. The EKI method has become popular recently, because of the fact that it does not require derivatives of the underlying model for optimization but at the same time enjoys provable convergence results. Applications have been so far, in particular, in oceanography~\cite{EvensenVanLeeuwen1996}, reservoir modeling~\cite{Aanonsen2009}, weather forecasting~\cite{McLaughlin2014}, milling process~\cite{SchwenzerViscontiEtAl2019}, process control~\cite{Teixeira2010}, geophysical applications~\cite{Keller2021,Muir2019,Tso2021}, physics~\cite{Li2020} and also machine learning~\cite{Haber2018NeverLB,KovachkiStuart2019,Alper}. The literature on Kalman filtering is very rich and we can not review this in detail here, but refer to the reference for further details. Our focus is on the reformulation of the EnKF for solving inverse problems as outlined below, in Section~\ref{ssec:intro:outline}.

\subsection{Formulation of the ensemble Kalman inversion} \label{ssec:enkf}

In order to present the mathematical formulation of the EKI method, we denote by $\G:X \to Y$  the given (nonlinear) forward operator between finite dimensional Hilbert spaces $X=\R^d$, $d\in\mathbb{N}$, and $Y=\R^K$, $K\in\mathbb{N}$. Consider the inverse problem or parameter identification problem of the type
\begin{equation} \label{eq:noisyProb}
\mbox{find } \ \vec{u}\in X \ \mbox{ s.t. } \ \G(\vec{u}) = \vec{y} + \vecsym{\eta}  \in Y.
\end{equation}
Throughout the paper $\vec{u}$ is referred to as the (unknown) control, whereas $\vec{y}$ represents the data measurements (that are perturbed by noise $\vecsym{\eta}$).
In applications one typically has $d\gg K$. The perturbations due to errors in the observations is modeled by $\vecsym{\eta}$ whose distribution is explicitly known. We assume that the noise is normally distributed with given covariance matrix $\vecsym{\Gamma}^{-1}\in\R^{K\times K}$, namely we write $\vecsym{\eta} \sim \mathcal{N}(\vec{0},\vecsym{\Gamma}^{-1})$.

In order to solve the inverse problem~\eqref{eq:noisyProb}, the EKI considers a number $J$ of particles or ensemble members whose state is determined by an iterative update. The ensemble members are modeled  as realizations of the control $\vec{u}\in\R^d$, in the following combined in $\vec{U}=\left\{\vec{u}^{j} \right\}_{j=1}^J$, with $\vec{u}^j\in\R^d$, $j=1,\dots,J$. The iteration index is denoted by $n$ and the collection of the ensemble members by $\vec{U}^n=\{\vec{u}^{j,n}\}_{j=1}^J$, $\forall\,n\geq 0$. 
%Then, 

%states are defined by evolving  of a possibly regularized minimization problem that compromises between the background estimate of the given model and additional information  $R$. 
%Note that $R\equiv 0$ is explicitly possible in the method. 
Then, at iteration $n+1$ the EKI update is given by 
%first order necessary optimality condition:
\begin{equation} \label{eq:updateEnKF}
\begin{aligned}
\vec{u}^{j,n+1} &= \vec{u}^{j,n} + \vec{C}_{\G}(\vec{U}^n) \left( \vec{D}_{\G}(\vec{U}^n) + \frac1{\Delta t} \vecsym{\Gamma}^{-1} \right)^{-1} (\vec{y} - \G(\vec{u}^{j,n}) ) \\
%{\vec{y}^{j,n+1}} &= \vec{y} + \vecsym{\xi}^{j,n+1}
\end{aligned}
\end{equation}
for each $j=1,\dots,J$, where $\Delta t\in\R^+$ is a parameter and where the ensemble update~\eqref{eq:updateEnKF} depends on covariance matrices:
\begin{equation} \label{eq:covariance}
\begin{aligned}
\vec{C}_{\G}(\vec{U}^n) &= \frac{1}{J} \sum_{k=1}^J \left(\vec{u}^{k,n}-\overline{\vec{u}}^n\right) \left(\G(\vec{u}^{k,n})-\overline{\G}^n\right)^\intercal \in \R^{d\times K} \\
\vec{D}_{\G}(\vec{U}^n) &= \frac{1}{J} \sum_{k=1}^J \left(\G(\vec{u}^{k,n})-\overline{\G}^n\right) \left(\G(\vec{u}^{k,n})-\overline{\G}^n\right)^\intercal \in \R^{K\times K}
\end{aligned}
\end{equation}
where we have denoted with $\overline{\vec{u}}^n$ and $\overline{\G}^n$ the mean of $\vec{U}^n$ and $\G(\vec{U}^n)$, respectively, namely
$$
\overline{\vec{u}}^n = \frac{1}{J} \sum_{j=1}^J \vec{u}^{j,n}, \quad \overline{\G}^n = \frac{1}{J} \sum_{j=1}^J \G(\vec{u}^{j,n}).
$$
Then, it can be proven~\cite{iglesiaslawstuart2013} that  
\begin{equation} \label{eq:originalMinProblem}
\lim_n  \overline{\vec{u}}^n = \arg\min_{\vec{u}} \frac12 \left\| \vecsym{\Gamma}^{\frac12} (\vec{y} - \G(\vec{u})) \right\|^2_Y. 
%+ R(\vec{u};\G(\vec{U}^n),\vec{U}^n).
\end{equation}
It is worth to mention that in the original formulation each observation or measurement is perturbed by additional additive noise at each iteration. 
The EKI satisfies the subspace property~\cite{iglesiaslawstuart2013}, i.e.,~the ensemble iterates stay in the subspace spanned by the initial ensemble. As consequence, the natural estimator for the solution of the inverse problem is provided by the mean of the ensemble. 
%The case $R=0$ arises in the continuous--time limit which has been recently proposed. This limit and its stabilization will be part of the following review. 

In recent years, the EKI was also studied as technique to solve inverse problems in a Bayesian framework. For instance see the works~\cite{ernstetal2015,Stuart2019MFEnKF} and the references therein. The analysis of the method is proven to have a comparable accuracy with traditional least--squares approaches to inverse problems~\cite{iglesiaslawstuart2013}. The method  approximates a specific Bayes linear estimators and it is able to provide  an approximation of the posterior measure. For a detailed discussion we refer  to~\cite{Apteetal2007,leglandmonbettran2009}. In this work, we keep the attention on the classical approach which aims to solve the inverse problem through an optimization point--of--view, see \eqref{eq:originalMinProblem}. Additional properties of the EKI method, continuous--time limits~\cite{schillingsetal2018,SchillingsPreprint,ChadaStuartTong2019,schillingsstuart2017,schillingsstuart2018}, i.e., $n\to \infty$  and mean--field limits on the number of the ensemble  members~\cite{CarrilloVaes,DingLi2019,Stuart2019MFEnKF,HertyVisconti2019}, i.e., $J\to \infty$ have been recently been developed and will be reviewed in more detail below.

\subsection{Structure of the paper} \label{ssec:intro:outline}

The remainder of this paper is organized as follows. In Section~\ref{sec:limitsEnKF} we review the continuous formulations of the EKI method which lead to a preconditioned gradient descent system and to a Vlasov--type partial differential equation. In Section~\ref{sec:multiobj} and in Section~\ref{sec:stabilization}, instead, we present two new formulations of the EKI method for multi--objective inverse problems and for globally asymptotically convergence to the target solution, respectively. Finally, we draw conclusions and perspectives in Section~\ref{sec:conclusion}. 

\section{Continuous limits of the ensemble Kalman inversion} \label{sec:limitsEnKF}

The continuous in time limit reduces the discrete update to a coupled system of ordinary differential equations. This limit has been performed in different recent publications, starting from~\cite{schillingsstuart2017} to more recent formulations, e.g.~see~\cite{Chada2020} for the hierarchical EKI. In particular, in~\cite{schillingsstuart2017} it has been shown that continuous in time limit results, in case of a linear forward model $\mathcal{G}$, to a gradient flow structure. This gradient flow provides a solution to the inverse problem~\eqref{eq:noisyProb} by minimizing the least--squares functional
\begin{equation} \label{eq:leastSqFnc}
\Phi(\vec{u},\vec{y}) := \frac12 \left\| \vecsym{\Gamma}^{\frac12} (\vec{y} - \G(\vec{u})) \right\|^2_Y.
\end{equation}
Observe, however, that in the continuous limit~\cite{schillingsstuart2017} the regularization term originally present in~\eqref{eq:originalMinProblem} vanishes for certain scalings. Although typically the analysis of the continuous in time EKI focuses on linear forward models, there are recent results on the EKI formulations in nonlinear settings~\cite{DingLiLu}.

\subsection{Continuous--time limit} \label{sec:limitsEnKF:time}
The continuous--time limit was firstly proposed in~\cite{schillingsstuart2017}: consider the parameter $\Delta t$ as an artificial time step for the discrete iteration, i.e.~$\Delta t \sim N_t^{-1}$ with $N_t$ being the maximum number of iterations and define  $\vec{U}^n \approx \vec{U}(n\Delta t)=\left\{\vec{u}^{j}(n\Delta t) \right\}_{j=1}^J$ for $n\geq 0$. Computing the limit $\Delta t\to 0^+$ one obtains
\begin{equation} \label{eq:continuousEnKF1}
\begin{aligned}
\frac{\mathrm{d}}{\mathrm{d}t} \vec{u}^j &= \vec{C}_{\G}(\vec{U}) \vecsym{\Gamma} \left( \vec{y} - \G(\vec{u}^j) \right), \quad j=1,\dots,J \\
\vec{C}_{\G}(\vec{U}) &= \frac{1}{J} \sum_{k=1}^J \left(\vec{u}^{k}-\overline{\vec{u}}\right) \left(\G(\vec{u}^{k})-\overline{\G}\right)^\intercal
\end{aligned}
\end{equation}
with initial condition $\vec{U}(0) = \vec{U}^0$. Note that within this limit the noise is scaled with $\frac{1}{\Delta t}$ which allows for the continuous time limit. Further, the term $\vec{D}_{\G}$ vanishes leading to possibly unstable dynamics~\cite{HertyVisconti2019,ArmbrusterHertyVisconti2022}. 

However, in the case of $\G$ linear, i.e.,~$\G(\vec{u})=\vec{G} \vec{u}$, with $\vec{G}\in\R^{K\times d}$, the~\eqref{eq:continuousEnKF1} can be reformulated in terms of the gradient  $\nabla \Phi$ as a  gradient flow:
\begin{equation} \label{eq:gradientEq}
\begin{aligned}
\frac{\mathrm{d}}{\mathrm{d}t} \vec{u}^j &= - \vec{C}(\vec{U}) \nabla_\vec{u} \Phi(\vec{u}^j,\vec{y}), \quad j=1,\dots,J \\
\vec{C}(\vec{U}) &= \frac{1}J \sum_{k=1}^J (\vec{u}^k-\overline{\vec{u}}) ( \vec{u}^k - \overline{\vec{u}} )^\intercal. 
\end{aligned}
\end{equation}
Since $\vec{C}(\vec{U})$ is positive semi--definite we obtain 
\begin{equation} \label{eq:decreasePhi}
\frac{\mathrm{d}}{\mathrm{d}t} \Phi(\vec{u}(t),\vec{y}) = \frac{\mathrm{d}}{\mathrm{d}t} \frac12 \left\| \vecsym{\Gamma}^{\frac12} \left(\vec{y}-\vec{G}\vec{u}\right) \right\|^2 \leq 0.
\end{equation}
Although the forward operator is assumed to be linear, the gradient flow is nonlinear. For further details and properties of the gradient descent equation~\eqref{eq:gradientEq} we refer to~\cite{schillingsstuart2017,schillingsstuart2018}.  In particular, we emphasize that the subspace property of the {EKI} also holds for the continuous dynamics and the following important result on the velocity of the collapse of the ensembles towards their mean in the large time limit, cf.~Theorem~3 in~\cite{schillingsstuart2017,SchillingsPreprint}:
$$
\left\|   \Gamma^{\frac12} G(\vec{u}^j(t)-\overline{\vec{u}}(t))\right\| = O(Jt^{-1}).
$$
%For  well--posedness of~\eqref{eq:gradientEq} and other properties we refer also to~\cite{SchillingsPreprint}.

\subsection{Mean--field limit} \label{sec:limitsEnKF:MFEnKF}

By definition, the {EKI} method considers a finite ensemble size $J<\infty.$ The behavior of the method in the limit of infinitely many ensembles can be studied via mean--field limit in analogy with the classical mean--field derivation of multi--agent systems~\cite{CarrilloFornasierToscaniVecil2010,Golse,jabin2014review, PareschiToscaniBOOK}. In the case of a linear foward model, the limit leads to a Vlasov--type gradient flow PDE. 
\begin{equation} \label{eq:kineticFromEnKF}
\partial_t f(t,\vec{u}) - \nabla_\vec{u} \cdot \left( \vec{C}(f) \nabla_\vec{u} \Phi(\vec{u},\vec{y}) f(t,\vec{u}) \right) = 0, \; f(t,0)=f_0(\vec{u})
\end{equation}
for a compactly supported on $\R^d$ probability density $f$ of $\vec{u}$ at time $t$ 
denoted by
\begin{equation} \label{eq:kineticf}
f = f(t,\vec{u}) : \R^+ \times \R^d \to \R^+.
\end{equation}
The initial probability density distribution is denoted by $f_0.$ The operator $\vec{C}(f)$ is the mean--field limit of the covariance of the ensemble and can be written in terms of moments of $f$ as
\begin{equation} \label{eq:covarianceMeanField}
\vec{C}(f) = \vec{E}(t) - \vec{m}(t) \vec{m}^\intercal(t) \geq 0,
\end{equation}
where $\vec{m}\in\R^{d}$ and $\vec{E}\in\R^{d\times d}$ are defined, respectively, as
\begin{equation} \label{eq:moments}
\vec{m}(t) = \int_{\R^d} \vec{u} f(t,\vec{u}) \mathrm{d}\vec{u}, \quad \vec{E}(t) = \int_{\R^d} \vec{u} \vec{u}^\intercal f(t,\vec{u}) \mathrm{d}\vec{u}.
\end{equation}
For the rigorous mean--field derivation and analysis of the EKI we refer to~\cite{CarrilloVaes,DingLi2019}. Equation~\eqref{eq:kineticFromEnKF} is a nonlinear transport equation arising from non--linear gradient flow interactions and in~\cite{CarrilloVaes,ArmbrusterHertyVisconti2022} it is observed that the counterpart of~\eqref{eq:decreasePhi} holds at the kinetic level. In fact, for 
$$
\mathcal{L}(f,\vec{y}) = \int_{\R^d} \Phi(\vec{u},\vec{y}) f(t,\vec{u}) \mathrm{d}\vec{u},
$$
we obtain 
$$
\frac{\mathrm{d}}{\mathrm{d}t} \mathcal{L}(f,\vec{y}) = \int_{\R^d} \Phi(\vec{u},\vec{y}) \partial_t f(t,\vec{u}) \mathrm{d}\vec{u} = - \int_{\R^d} (\nabla_\vec{u} \Phi(\vec{\vec{u}},\vec{y}) )^T \vec{C}(f) \nabla_\vec{u} \Phi(\vec{\vec{u}},\vec{y}) \mathrm{d}\vec{u} \leq 0,
$$
since $\vec{C}(f)$ is positive semi-definite. In particular, $\mathcal{L}(f,\vec{y})$ is  strictly decreasing unless $f$ is a Dirac measure.  Also, for $f(\vec{u})=\delta(\vec{u}-\vec{u}^*)$ for  $\vec{u}^*=argmin_{\vec{u}\in\R^d}\Phi(\vec{u},\vec{y})$ provides a steady solution of the continuous--limit formulation, but the converse is not necessarily true.  In fact, all Dirac distributions, satisfy $\vec{C}(f)=0$ and hence provide steady solutions of~\eqref{eq:kineticFromEnKF}. 
Convergence to the  distribution $f(\vec{u})=\delta(\vec{u}-\vec{u}^*)$  has been proven to be linear in~\cite{CarrilloVaes}: $\left\| \vec{C}(f) \right\| = O(t^{-1})$.

The mean--field interpretation of the EKI has allowed to design computationally efficient methods based on the mean--field formulation~\cite{AlbiPareschi2013,HertyVisconti2019}. In particular, it is possible to use a large number of particles which guarantees significantly better reconstructions of the unknown control, cf.~Section 5 in~\cite{HertyVisconti2019}.

\section{Multi--objective ensemble Kalman inversion} \label{sec:multiobj}

The EKI can also be extended to treat also multi-objective optimization problems within a weighted function approach. 
Here, a vector of controls has to be determined for competitive models $\mathcal{G}_i$ for $i=1,\dots,l$ and  given observational data: 
\begin{equation}\label{eq:coupled_ip}
\vec{y}_i=\G_i(\vec{u}) + \vecsym{\eta}_i \quad i=1,\dots,l
\end{equation}
for $l$ models $\G_i:X \to Y$ and $l$ observations $\vec{y}_1,\dots,\vec{y}_l \in Y$, where $ \vecsym{\eta}_i$ is  observational noise. A solution to  \eqref{eq:coupled_ip} can be obtained e.g. using a multi--objective optimization  \cite{ehrgott2005multicriteria, miettinen2012nonlinear,pardalos2017non}:

{\begin{equation}\label{eq:enkfo2}
	\min_{\vec{u} \in X}  \left( \| \Gamma^{\half} \left( \vec{y}_1-\G_1(\vec{u}) \right) \|, \dots,  \| \Gamma^{\half} 
	\left( \vec{y}_l-\G_l(\vec{u}) \right) \|
	\right).
	\end{equation}}

Solution in this framework is related to the notion of  Pareto optimality \cite{pardalos2017non} that defines a concept of minimum for the vector--valued optimization problem \eqref{eq:enkfo2}.

\begin{definition}\label{def1}
	A point $\vec{u}^*\in \R^d$ is called Pareto optimal if and only if there exists no point $\vec{u} \in \R^d$ such that $\G_i(\vec{u}) \le \G_i(\vec{u}^*)$ for all $i=1,2,\dots,l$ and $\G_j(\vec{u})\le \G_j(\vec{u}^*)$ for at least one $j\in \{1,2,\dots,l\}$.
\end{definition} 
The set $\mathcal{S}_U$ of all  $\vec{u}^*$ fulfilling Definition~\ref{def1} is called Pareto set, while its representation in the space of objectives  $\mathcal{S}_G := \{ \left(\vec{y}_i-\G_i(\vec{u})\right)_{i=1}^l : \vec{u} \in \mathcal{S} \}$ is called Pareto front. An approximation of $\mathcal{S}_G$ can be recovered following an approach based on weighted function method \cite{miettinen2012nonlinear}. Let   $\mathbf{1}=(1,\dots,1)^T$ and let $\mathbf{\mathbb{\lambda}} \in   \Lambda$ be a fixed vector in the set 
\begin{equation}\label{lambda} \Lambda:=\{ \mathbf{\mathbb{\lambda}} \in \R^l_+: \mathbf{\mathbb{\lambda}} \cdot \mathbf{1}=1   \}. 
\end{equation}
Define the weighted objective function and the weighted observations as 
\begin{equation}\label{eq:g}
\G(\vec{u},\mathbf{\mathbb{\lambda}}):= \sum\limits_{i=1}^l \mathbf{\mathbb{\lambda}}_{i}  \G_i(\vec{u}): X\times \Lambda \to Y \ \mbox{ and } \ \vec{y}= \sum\limits_{i=1}^l \mathbf{\mathbb{\lambda}}_{i}  \vec{y}_i.
\end{equation}
An approximation to the Pareto front $\mathcal{S}_U$ is then obtained by 
%\begin{eqnarray}\label{paretof}
${P}:=\{ \vec{u}^*(\mathbb{\lambda}) : \mathbb{\lambda} \in \Lambda  \}$, %\label{pareto1}
%\end{eqnarray}
where for each $\lambda$ 
\begin{equation}\label{opt1}
\vec{u}^*(\mathbb{\lambda}) = \arg\min_{\vec{u} \in X} \Phi(\vec{u},\mathbb{\lambda}), \quad \Phi(\vec{u},\vec{y},\mathbb{\lambda})=\frac{1}{2} \norm{\Gamma ^{\half} \sum_{i=1}^l\mathbb{\lambda}_i \left( \vec{y}_i-\G_i(\vec{u}) \right)}^2 \; \forall \lambda \in \Lambda.
\end{equation}
In case of a convex problem, $S_U = P$, {see \cite[Theorem 3.1.4]{miettinen2012nonlinear}. In theory the previous problem \eqref{opt1} needs to be solved for all $\lambda\in\Lambda.$ 
	
	Using a mean--field approach as in the last section, allows for an analysis on the dependence of $u^*(\lambda)$ on $\lambda$ which in turn is used as  adaptive grid  on $\Lambda$.  The evolution of the formal sensitivity of the mean-field description $f$ of the particle distribution with respect to $\lambda_i$ is given by
	\begin{equation} \label{eq:der}
	\begin{aligned}
	0 & =  \partial_t \partial_{\lambda_i} f(t,\vec{u},\lambda) - \nabla_\vec{u} \Big( \partial_{\lambda_i} \vec{C}(f) \nabla_\vec{u} \Phi(\vec{u},\vec{y},\lambda) f(t,\vec{u},\lambda) + \\ & \qquad \vec{C}(f)\ \partial_{\lambda_i} (\nabla_\vec{u}  \Phi(\vec{u},\vec{y},\lambda)) f(t,\vec{u},\lambda) + \vec{C}(f)\ \nabla_\vec{u}  \Phi(\vec{u},\vec{y},\lambda) \partial_{\lambda_i} f(t,\vec{u},\lambda) \Big),
	\end{aligned}
	\end{equation}
	for zero initial data. The set of equations \eqref{eq:der} for $i=1,\dots,l$ is defined on the extended phase space 
	$\mathbb{R} \times X\times \Lambda$ and therefore computationally infeasible. However, the Pareto set is given as moment of $f$ where the first moment $\vec{m}$ depends additionally on $\lambda:$
	\begin{equation}\label{eq:par}
	{P}(t)=\left\{ \vec{m} (\lambda,t):= \int_X \vec{u} \; \mathrm{d} f(t,\vec{u},\mathbb{\lambda}): \;   \mathbb{\lambda} \in \Lambda  \right\}.
	\end{equation}
	Similarly, moments of \eqref{eq:der} can be defined leading to a set of ordinary differential equations for $\nabla_\lambda \vec{m} (\lambda,t)$  \cite[Lemma 2.3]{hertyIac}. This in turn allows to define an adaptive grid $\{ \lambda^k\, k=1,\dots \}  \subset \Lambda$: Let for a fixed $\overline{\lambda}=\lambda^{k-1}$ the  corresponding optimal parameter be approximated  by $\vec{m}(\overline{\lambda},T)$ for some $T$ fixed and sufficiently large. Then,  consider the following Taylor expansion 
	\begin{align}\label{formula}
	\vec{m}(\overline{\lambda}+\Delta \lambda,T) \approx  \vec{m}(\overline{\lambda},T) + \Delta \lambda \cdot \nabla \vec{m}(\overline{\lambda},T), \mbox{ and } \Delta \lambda = \lambda^{k+1}-\overline{\lambda}. 
	\end{align}
	Reformulating \eqref{formula} allows to obtain $\lambda^{k+1} \in \Lambda$ adaptively based on the  approximation on the Pareto set $P(t)$. It also yields an  estimate on the norm of the update $\Delta \lambda$ on  an approximation of $S_U$ with given  tolerance $\delta>0$ by  $\| \Delta \lambda \| \| \nabla \vec{m}(\lambda^k ) \| \leq \delta.$

	\subsection{Numerical experiment} \label{sec:multiobj:numerics}
	
	In the numerical experiment we show that the adaptive strategy leads to results that  approximate the Pareto front $S_G$ very well with only a few discretization points $\lambda^k, k=1,\dots,K$.  We set $l=2$ so that $\Lambda$ is parameterized by a single parameter $\lambda \in [0,1]$, i.e. $\G=\lambda \G_1+(1-\lambda)\G_2$. 
	Then, we consider two non convex  functions $\G_1,\G_2: \R^2\to \R^2$ as in \cite{deb2002fast}  
	\begin{equation}
	\G_1(u_1,u_2)= 1-e^{-\left(u_1-\frac{1}{\sqrt{2}}\right)^2-\left(u_2-\frac{1}{\sqrt{2}}\right)^2}, \quad \G_2(u_1,u_2)= 1-e^{-\left(u_1+\frac{1}{\sqrt{2}}\right)^2-\left(u_2+\frac{1}{\sqrt{2}}\right)^2},
	\end{equation} 
	and  $\vec{y}_i=0$ and $ \vecsym{\eta}_i=0$, for $i=1,2$.  As further parameters we use $J=25$ particles sampled from the uniform distribution $U_0\sim \mathcal{U}([-2,2]^2)$,  the tolerance is set $\delta=5\cdot 10^{-3}$, $T_{fin}=10$, $\Gamma=\mathbbm{1}$ and $K=22$.
	
	Even so, the theoretical results have been proven in the linear case \cite[Sec. 3]{hertyIac}, they are applied here in a nonlinear framework. We compare a naive choice for the discretization of $\Lambda$ using an equidistant grid (direct approach) with the outlined adaptive strategy.
	
	We observe that the solution obtained with the adaptive approach covers a larger part of the Pareto front, showing additionally a relatively sharper resolution  compared with the direct approach, see Figure~\ref{fig:paretofront}.
	
	\begin{figure}[t!]
		\centering
		\includegraphics[width=\textwidth]{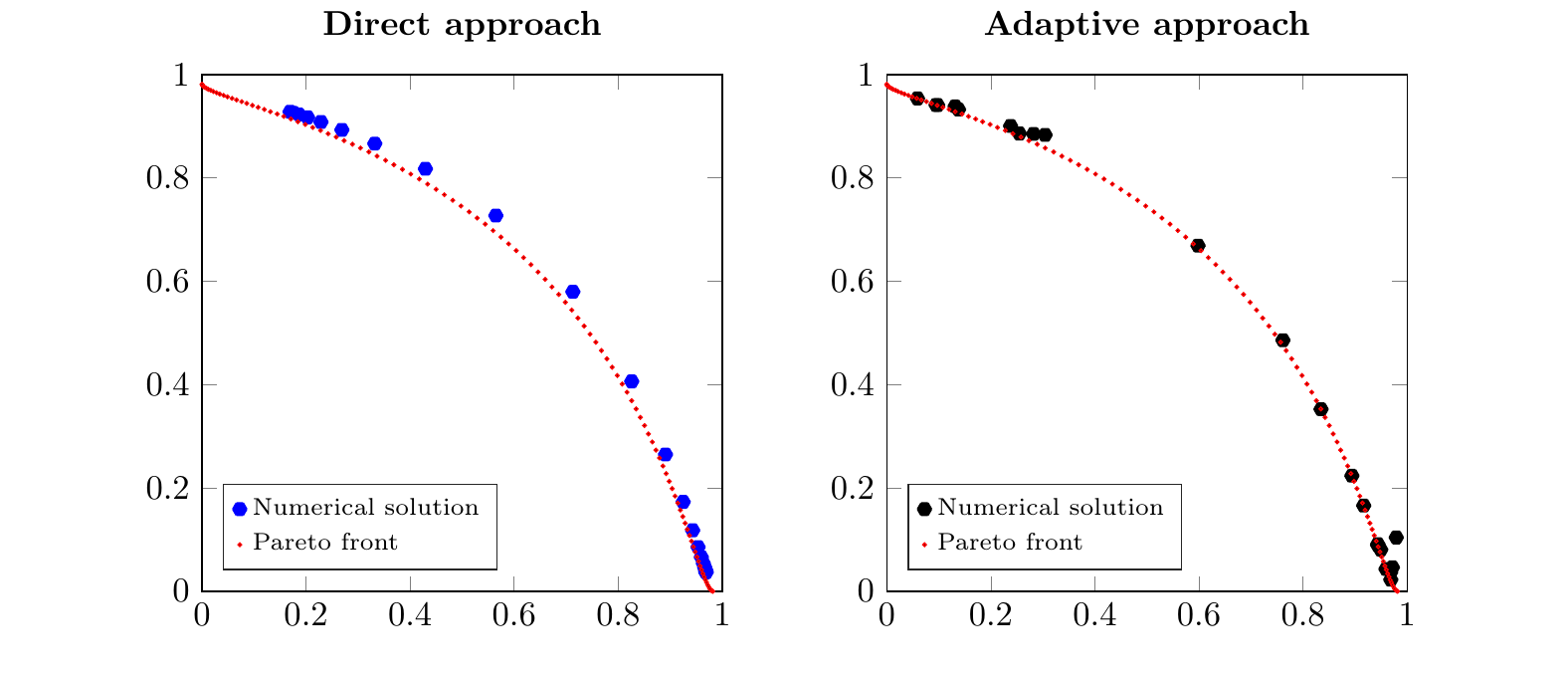}
		\caption{Approximation of the Pareto front  with the direct approach (left) and the adaptive approach (right).}\label{fig:paretofront}
	\end{figure}

	%In order to highlight the difference between the two strategies, 
	Moreover, the approximation of the Pareto set in Figure~\ref{fig:paretoset} shows the expected behavior. Here, the  adaptive strategy yields a cloud of points relatively close to the (analytically known) Pareto set $S_U$ compared with the direct approach. 
	\begin{figure}[]
		\centering
		\includegraphics[width=0.49\textwidth]{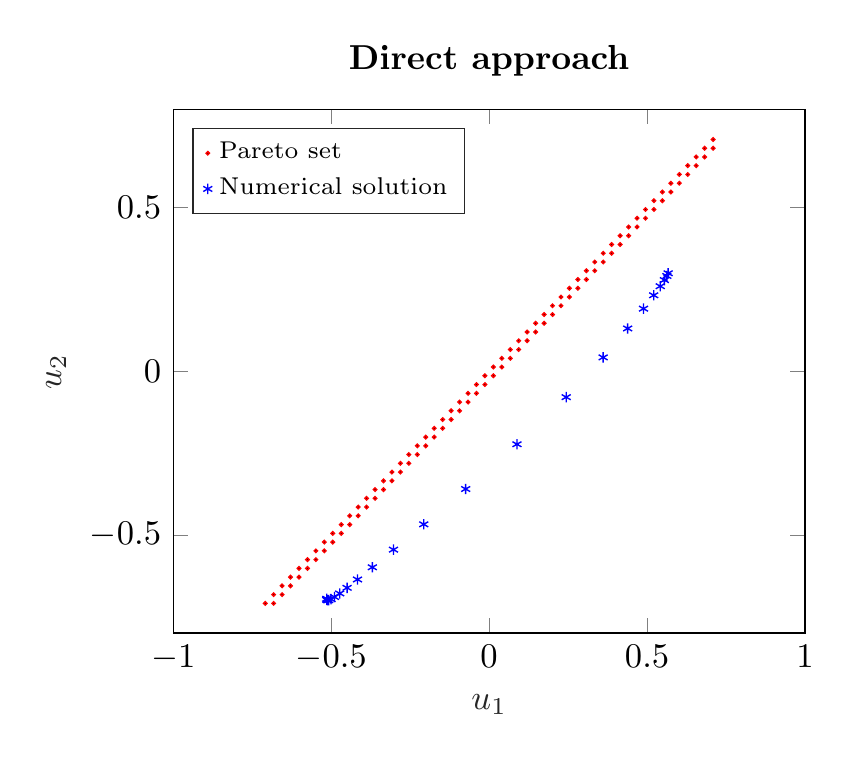}
		\includegraphics[width=0.49\textwidth]{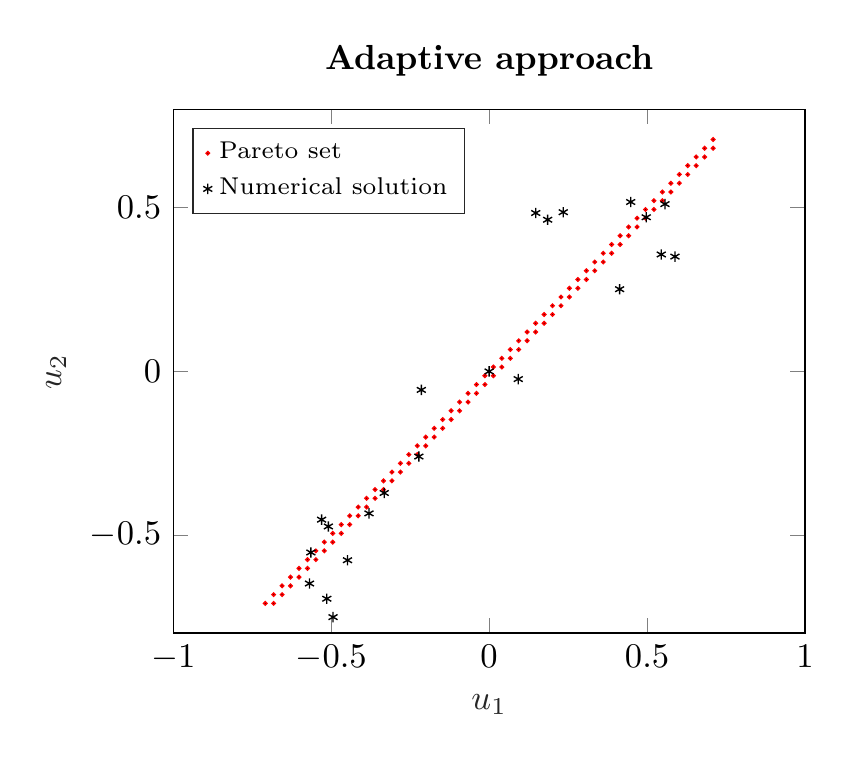}
		\caption{Approximation of the Pareto set with the direct approach (left) and the adaptive approach (right).}\label{fig:paretoset}
	\end{figure}
	%One may notice that in the adaptive framework the cloud of points is relatively close to the Pareto set, indeed several points lay on it, while the approximation obtained with the direct approach is quite far.
	
	\section{Stabilized continuous limit of the ensemble Kalman inversion} \label{sec:stabilization}
	
	In the continuous--time limit the  term  $\vec{D}_{\mathcal{G}}$ present in the discrete formulation vanish due to scaling.  This consideration inspired~\cite{ArmbrusterHertyVisconti2022}, where a  stability analysis of the moment system of the time--continuous EKI~\eqref{eq:gradientEq} is performed. Therein, it has been established that the system has infinitely many non--hyperbolic Bogdanov--Takens equilibria leading to several undesirable consequences. The latter are structurally unstable, i.e.,  sensitive to small  perturbations. Since  those equilibria lie on the set where the preconditioner $\vec{C}$ collapse to zero, low order of convergence in time holds true. Further, numerical  approximations  may  push the trajectory in the unfeasible region of the phase space or get the method stuck in the wrong equilibrium. 
	
	These considerations led to a modified formulation of the method  is globally asymptotically stable by introducing the regularization term $R$ to the dynamics.  More precisely, given $\vecsym{\Sigma}\in\R^{d\times d}$ symmetric and full--rank, in particular positive definite, in~\cite{ArmbrusterHertyVisconti2022} it is proposed to consider the following general discrete dynamics for each ensemble member $j=1,\dots,J$ in the case of a linear model:
	\begin{equation} \label{eq:stableLinearEnKF}
	\begin{aligned}
	\frac{\mathrm{d}}{\mathrm{d}t} \vec{u}^j &= - \tilde{\vec{C}}(\vec{U}) \nabla_\vec{u} \Phi(\vec{u}^j,\vec{y}) + R(\vec{U}), \\
	R(\vec{U}) &= \beta \tilde{\vec{C}}(\vec{U})(\vec{u}^j-\bar{\vec{u}}), \quad 
	\tilde{\vec{C}}(\vec{U}) = \vec{C}(\vec{U}) + (1-\alpha)\vecsym{\Sigma},
	\end{aligned}
	\end{equation}
	with parameters $\alpha,\beta\in\R$. The choices $\alpha=1$ and $\beta=0$ yield the continuous--time limit~\eqref{eq:gradientEq} for the original EKI. The modified dynamics~\eqref{eq:stableLinearEnKF} differs from~\eqref{eq:gradientEq} in the formulation of the preconditioner $\tilde{\vec{C}}(\vec{U})$ and in the presence of the additive term $R(\vec{U})$. The new preconditioner is related to an inflation of the covariance $\vec{C}(\vec{U})$ for $\alpha<0$. This modification allows  to stabilize the phase space of the moments.  The term $R(\vec{U})$, instead, has been shown to be an acceleration term for the convergence towards equilibrium. The  modified dynamical system~\eqref{eq:stableLinearEnKF} has also a mean--field interpretation:
	\begin{equation} \label{eq:stableMeanField}
	\begin{aligned}
	\partial_t f(t,\vec{u}) - \nabla_\vec{u} \cdot \left( \tilde{\vec{C}}(f) \left( \nabla_\vec{u} \Phi(\vec{u},\vec{y}) - \beta (\vec{u}-\vec{m}) \right) f(t,\vec{u}) \right)= 0,
	\end{aligned}
	\end{equation}
	where $\tilde{\vec{C}}(f)$ is the mean--field of $\tilde{\vec{C}}(\vec{U})$ leading to
	$
	\tilde{\vec{C}}(f) = \vec{E}(t) - \vec{m}(t)\vec{m}^\intercal(t) +(1-\alpha)\vec{\Sigma}.$
	
	The stability analysis of the moment equations is performed in the simplified setting where $K=d$ and $\vec{\Gamma}$, $\vec{G}$ are identity matrices. The $d+d^2$ dynamical system of the moments of~\eqref{eq:stableMeanField} is then
	\begin{equation*}
	\begin{aligned}
	\frac{\mathrm{d}}{\mathrm{d}t} \vec{m}(t) &= \vec{C} (\vec{y}-\vec{m}) + (1-\alpha) \vecsym{\Sigma} (\vec{y}-\vec{m})\\
	\frac{\mathrm{d}}{\mathrm{d}t} \vec{C}(t) &= -2\vec{C} \vec{C} - (1-\alpha) \vecsym{\Sigma} \vec{C} - (1-\alpha) \vec{C} \vecsym{\Sigma},
	\end{aligned}
	\end{equation*}
	and its  linearization $(\delta\vec{m},\delta\vec{E})$ at  target equilibrium $(\vec{m}^*,\vec{C}^*)=(\vec{y},\vec{0})$ fulfills
	\begin{align*}
	\frac{\mathrm{d}}{\mathrm{d}t} \delta\vec{m}(t) &=  - (1-\alpha) \vecsym{\Sigma}  \delta \vec{m}, \quad
	\frac{\mathrm{d}}{\mathrm{d}t} \delta\vec{C}(t) = - 4 (1-\alpha) \vecsym{\Sigma} \delta \vec{C}.
	\end{align*}
	For  $\alpha<1$ and $\vecsym{\Sigma}$ positive definite the target equilibrium is hyperbolic. This formal presentation of the role of the parameters is mathematically rigorous in~\cite{ArmbrusterHertyVisconti2022}, where also exponentially fast convergence to the target equilibrium is proven.
	
	\subsection{Numerical experiments} \label{sec:stabilization:numerics}
	
	We consider the inverse problem of identifying  the force function $u(x)$ of a linear elliptic equation in one spatial dimension assuming that noisy observation of the solution to the problem are available, e.g.~see~\cite{iglesiaslawstuart2013,schillingsstuart2017,HertyVisconti2019}.
	
	The problem is prescribed by the following one dimensional elliptic PDE for $p$ 
	\begin{equation} \label{eq:ellipticEq}
	-\frac{\mathrm{d}^2}{\mathrm{d}x^2} p(x) + p(x) = u(x), \quad x\in[0,\pi]
	\end{equation}
	subject to boundary conditions $p(0) = p(\pi) = 0$. To obtain measurement data we use the continuous  control $u(x)=\sin(8x)$. The problem is discretized using  a uniform mesh with  $d=K=2^8$ equidistant points $\{ x_i \}$ on the interval $[0,\pi]$. Denote by  $\vec{u}^\dagger_i = u(x_i)$ the evaluations of the control function $u(x)$ on the mesh $i=0,1\dots,d.$ Noisy observations $\vec{y}\in\R^K$ are obtained by
	$$
	\vec{y} = \vec{p} + \vecsym{\eta} = \vec{G} \vec{u}^\dagger + \vecsym{\eta},
	$$
	where $\vec{G}\in\R^{K\times d}$ is a first--order finite difference discretization of the PDE~\eqref{eq:ellipticEq}. For simplicity we assume that $\vecsym{\eta}$ is a Gaussian white noise, i.e., $\vecsym{\eta}\sim \mathcal{N}(0,\gamma^2 \vec{I})$ with $\gamma \in \R^+$ and $\vec{I} \in \R^{d\times d}$ is the identity matrix. We are interested in recovering an approximation to the  discrete control $\vec{u}^\dagger \in \R^d$ from the noisy observations $\vec{y}\in\R^K$.
	
	In Figure~\ref{fig:stabilization:solControl} we show the solution to this problem provided by the time--continuous limit of the original EKI and the stabilized formulation proposed in~\cite{ArmbrusterHertyVisconti2022}. Both method use $J=20$ ensemble members and a noise level of $\gamma=0.01$.  We observe that the stable EKI produces a qualitatively improved reconstruction of the control and  observation compared to the classical EKI. Moreover, as  expected by the analysis,  we observe that the stable EKI converges faster than the classical method, see Figure~\ref{fig:stabilization:analysis}. 
	%The simulation stops at the time where the stable EKI meets the discrepancy principle.
	
	\begin{figure}[t!]
		\centering
		\includegraphics[width=\textwidth]{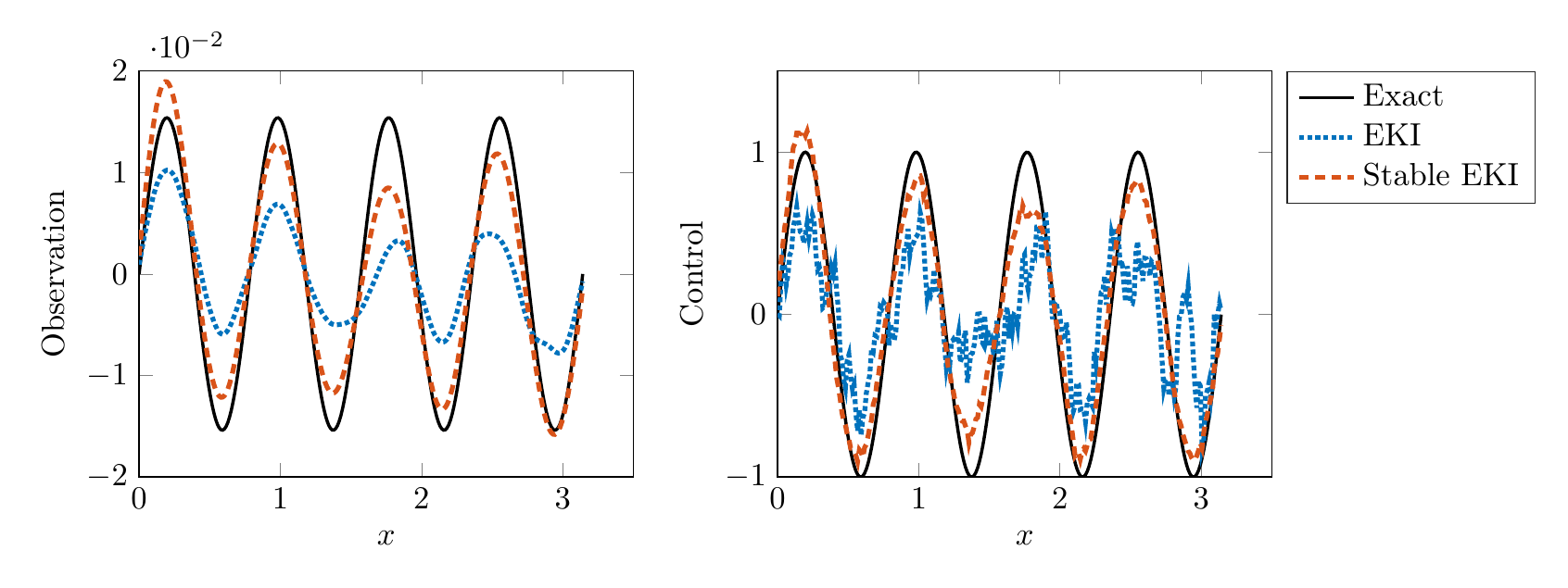}
		\caption{Reconstruction of the observation (left) and reconstruction of the control (right) for the continuous--time limit of the original EKI method and for the stabilized formulation}.\label{fig:stabilization:solControl}
	\end{figure}
	
	\begin{figure}[t!]
		\centering
		\includegraphics[width=\textwidth]{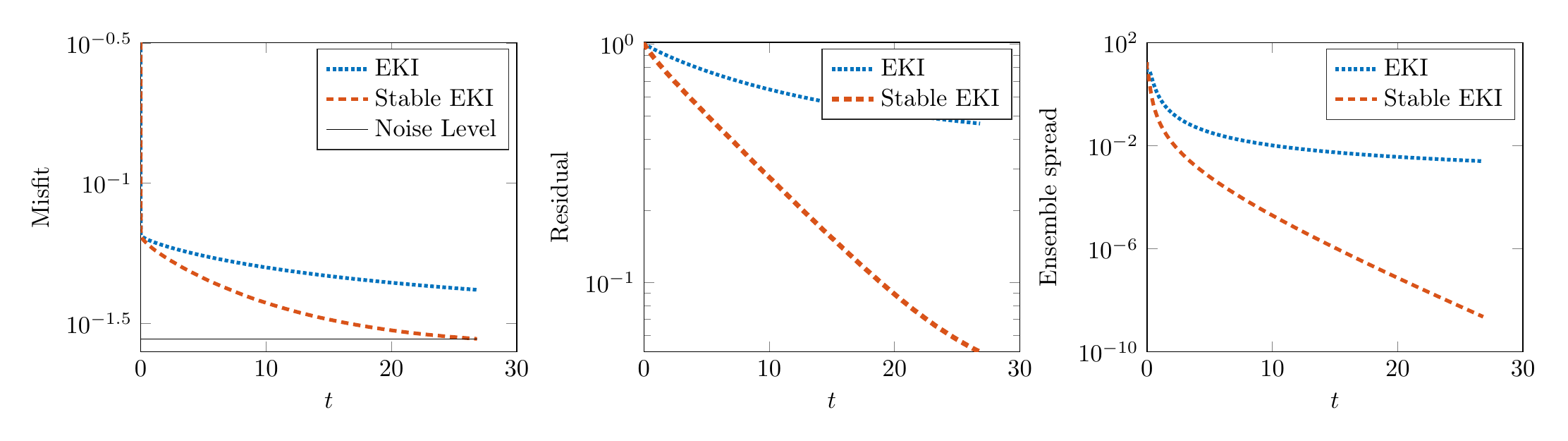}
		\caption{Behavior of the misfit, of the residual and of the ensemble spread around the mean in time, both for the continuous--time limit of the original EKI method and for the stabilized formulation in~\cite{ArmbrusterHertyVisconti2022}.\label{fig:stabilization:analysis}}
	\end{figure}
	
	\section{Conclusions} \label{sec:conclusion}
	
	An overview on the EKI and its current developments has been provided. The analytical properties have been investigated and, in particular,  the mean--field equation and its corresponding moment system has been presented.
	Two recent extensions of the EKI has been shown and discussed towards  coupled inverse problems and towards numerically stable formulations.  Further developments may involve a mixture between the two novelties presented and since many physical problems are subject to additional parameteric uncertainty, a suitable treatment of the then stochastic EKI might be of further interest. In case of large parameter spaces $X=\R^d$ with $d\gg1$ computational issues need to be addressed as well, since e.g. $\vec{C}_{\mathcal{G}}$ grows quadratic in $d.$ Furthermore, the outlined approach of time-continuous and mean-field limit is applicable to wider range of particle methods and might serve as a starting point for future investigation into nonlinear filtering from a mathematical perspective.
	
	\section*{Acknowledgments}
	
	The authors thank the Deutsche Forschungsgemeinschaft (DFG, German Research Foundation) for the financial support through 20021702/GRK2326, 333849990/IRTG-2379,
	HE5386/15,18-1,19-1,22-1,23-1 and under Germany’s Excellence Strategy EXC-2023 Internet of Production 390621612. The funding through HIDSS-004 is acknowledged.

	G.V.~is member of the ``National Group for Scientific Computation (GNCS-INDAM)'' and acknowledges support by MUR (Ministry of University and Research) PRIN2017 project number 2017KKJP4X.

\bibliography{references}  
\bibliographystyle{abbrv}

\end{document}